\documentclass[11pt]{article}
\usepackage{a4wide}

\usepackage[style=numeric,sorting=none]{biblatex}
\usepackage[american]{babel}
\usepackage[fleqn]{amsmath}
\usepackage{amssymb}
\usepackage{hyperref}
\usepackage{tabularx}
\usepackage{ltablex}
\usepackage{csquotes}
\usepackage{graphicx}
\usepackage{placeins}
\usepackage{eurosym}
\usepackage{comment}
\keepXColumns
\usepackage[dvipsnames]{xcolor}
\title{Tulipa Energy Model: Mathematical Formulation}
\author{Diego A. Tejada-Arango$^a$, Germán Morales-España$^a$, Lauren Clisby$^a$, Ni Wang$^a$, \\ 
Abel S. Siqueira$^b$, Ali Subayu $^b$, Laurent Soucasse $^b$, and Zhi Gao $^c$ \\ \\
$^a$ \textit{TNO - Energy Transition Studies} \\
$^b$ \textit{eScience Center} \\
$^c$ \textit{Utrecht University}
}
\date{\today}

\begin{document}

\maketitle

\section*{Abstract}
\textit{Tulipa aims to optimise the investment and operation of the electricity market, considering its coupling with other sectors, such as hydrogen and heat, that can also be electrified. The problem is analysed from the perspective of a central planner who determines the expansion plan that is most beneficial for the system as a whole, either by maximising social welfare or by minimising total costs. The formulation provides a general description of the objective function and constraints in the optimisation model based on the concept of energy assets representing any element in the model. The model uses subsets and specific methods to determine the constraints that apply to a particular technology or network, allowing more flexibility in the code to consider new technologies and constraints with different levels of detail in the future.}

\section{Tulipa Energy Model: An electricity market and sector-coupling model for investment and operation}
\begin{figure}
    \centering
    \includegraphics[width=\textwidth]{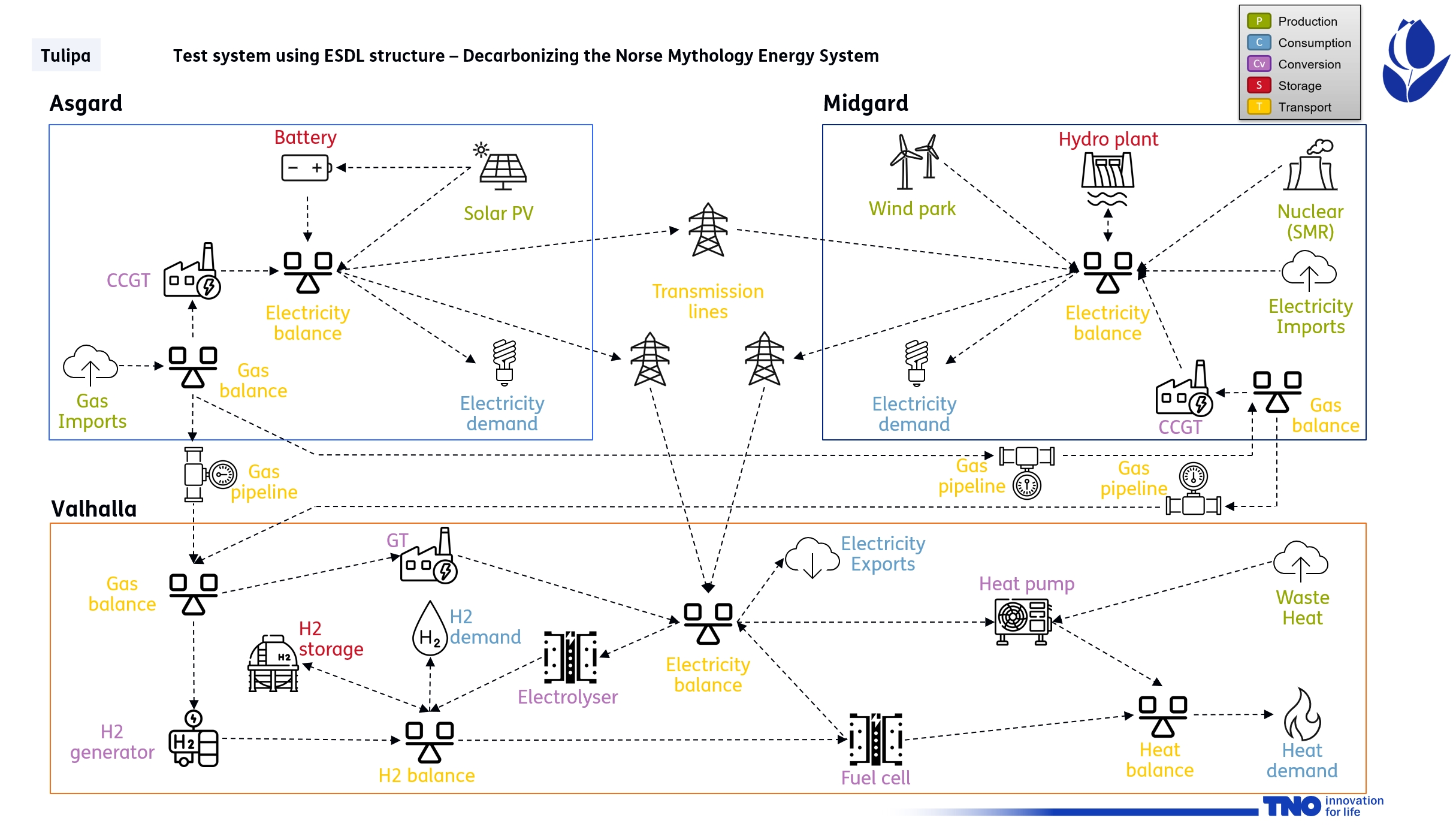}
    \caption{Test system in Tulipa Energy Model}
    \label{fig:test_system}
\end{figure}

In general terms, this document presents a classical problem of joint generation and transmission expansion planning in an electric energy system \cite{Conejo2016}. However, Tulipa's formulation is more generalised and includes other sectors coupled with the electricity sector. In addition, the formulation is based on Graph Theory \cite{West2000}, which provides a more flexible framework to model energy assets in the system as \textit{vertices} and flows between energy assets as \textit{edges}. Figure (\ref{fig:test_system}) shows an example of an energy system network coupling electricity, methane (i.e., gas), and hydrogen ($H_2$).

The following sections describe the notation, constraints, and objective function in the formulation used in the model.

\section{Notation}

\subsection{Sets for assets}
The formulation in Tulipa relies on the definitions of energy assets (i.e., vertices of the network), which have been taken from the Energy System Description Language (ESDL) \cite{ESDL}.

\begin{tabularx}{\textwidth}{| l | X | r |} 
\hline 
\textbf{Name} & \textbf{Description} & \textbf{Elements}\\ 
\hline 
\endhead 
$\mathcal{A}$      & Energy assets      & $a                \in     \mathcal{A}$ \\
$\mathcal{A}_{cv}$ & Conversion assets  & $\mathcal{A}_{cv} \subset \mathcal{A}$ \\
$\mathcal{A}_{p} $ & Production assets  & $\mathcal{A}_{p}  \subset \mathcal{A}$ \\
$\mathcal{A}_{c} $ & Consumption assets & $\mathcal{A}_{c}  \subset \mathcal{A}$ \\
$\mathcal{A}_{t} $ & Transport assets   & $\mathcal{A}_{t}  \subset \mathcal{A}$ \\
$\mathcal{A}_{s} $ & Storage assets     & $\mathcal{A}_{s}  \subset \mathcal{A}$ \\
\hline 
\end{tabularx}

As shown in the previous Table, energy assets, denoted by $a$, can belong to different subsets based on the definitions, such as conversion, production, consumption, transport, and storage assets. Here there are some examples:

\begin{itemize}
    \item Conversion ($\mathcal{A}_{cv}$): Power plants, fuel cells, heat pumps, etc.
    \item Production ($\mathcal{A}_p$): Solar PV panels, wind farms, etc.
    \item Consumption ($\mathcal{A}_c$): Electricity demand, heat consumption of a city, etc.
    \item Transport ($\mathcal{A}_t$): Electricity nodes and lines, gas networks, district heating, etc.
    \item Storage ($\mathcal{A}_s$): Battery, pumped storage, heat buffer, hydrogen storage, etc.

\end{itemize}

\begin{figure}[ht]
    \centering
    \href{https://energytransition.gitbook.io/esdl/esdl-concepts/design-principles}
        {\includegraphics[width=0.5\textwidth]{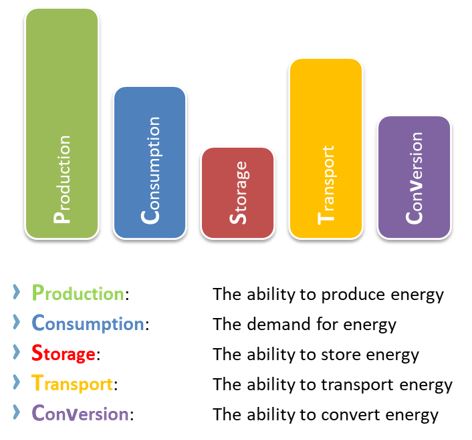}}        
    \caption{Energy assets in ESDL}
    \label{fig:esdl_energy_assets}
\end{figure}

General constraints can be created at the energy asset level ($a$), such as maximum capacity, and specific constraints can be applied to these subsets or groups, such as unit commitment constraints for the conversion units ($\mathcal{A}_{cv}$). This structure allows the model to be flexible in representation and levels of detail.

\subsection{Extra sets}
The formulation also uses more set definitions, such as:
\begin{tabularx}{\textwidth}{| l | X | r |} 
\hline 
\textbf{Name} & \textbf{Description} & \textbf{Elements}\\ 
\hline 
\endhead 
$\mathcal{A}_b$       & Assets with balance constraint method                             & $\mathcal{A}_b             \subseteq \mathcal{A}$ \\
$\mathcal{A}^{in}(a)$ & Assets that are inputs for asset $a$                              & $\mathcal{A}^{in}(a)       \subseteq \mathcal{A}$ \\
$\mathcal{A}^{out}(a)$& Assets that are outputs for asset $a$                             & $\mathcal{A}^{out}(a)      \subseteq \mathcal{A}$ \\
$\mathcal{Y}$         & Years                                                             & $y                         \in       \mathcal{Y}$ \\
$\mathcal{Y}_{m} $    & Milestone years                                                   & $\mathcal{Y}_{m}           \subseteq \mathcal{Y}$ \\
$\mathcal{K}     $    & Representative periods                                            & $k                         \in       \mathcal{K}$ \\
$\mathcal{T}     $    & Time steps for the lowest time resolution in the system           & $t                         \in       \mathcal{T}$ \\
$\tau(a)         $    & Time steps defined for each asset $a$                             & $\tau(a)                   \subseteq \mathcal{T}$ \\
$\mathcal{P}     $    & Auxiliary set for constraints with several time steps definitions & $p \in \mathcal{P}         \subseteq \mathcal{T}$ \\
$\mathcal{R}     $    & Types of system reserves                                          & $r                         \in       \mathcal{R}$ \\
$\mathcal{R}_{a}^{+} $& Upwards reserves provided by asset $a$                            & $\mathcal{R}_{a}^{+}       \subseteq \mathcal{R}$ \\
$\mathcal{R}_{a}^{-} $& Downwards reserves provided by asset $a$                          & $\mathcal{R}_{a}^{-}       \subseteq \mathcal{R}$ \\
$\Omega_{ap\tau} $    & Set to map between two temporal definitions $p$ and $\tau$        &                                                   \\
\hline 
\end{tabularx}

The model uses sets to define methods, which are used to determine constraints and detail levels for energy assets. For instance, depending on whether it belongs to $\mathcal{R}_{a}^{+}$ or $\mathcal{R}_{a}^{-}$, an energy asset may have a reserve requirement method. The balance method determines if the asset uses the balance constraint when it belongs to $\mathcal{A}_b$. The concept of methods has been previously defined in \cite{Helistö2019} and \cite{Ihlemann2022}.

\subsection{Parameters}
The model's parameters are listed in the table below. Unit measures and scales may vary based on input data.

\begin{tabularx}{\textwidth}{| l | X | l |} 
\hline 
\textbf{Name} & \textbf{Description} & \textbf{Units}\\ 
\hline 
\endhead 
$\overline{F}_{a}$         & Maximum asset capacity                                        & [MW]        \\
$\overline{F}^{max}_{ya}$  & Maximum potential capacity                                    & [MW]        \\
$\tilde{F}^{max}_{yakt}$   & Maximum energy asset production profile                       & [p.u.]      \\
$\tilde{F}^{min}_{yakt}$   & Minimum energy asset production profile                       & [p.u.]      \\
$\overline{U}_{ya}$        & Initial installed assets                                      & [-]         \\
$S_{yakt}^I$               & Energy asset inflows per representative                       & [MWh]       \\
$\eta_{a \alpha}$          & Efficiency when converting energy from $a$ to $\alpha$        & [p.u.]      \\
$R_{yrkt}$                 & Reserve requirement                                           & [MW]        \\
$D_{yakt}$                 & Demand of asset $a$                                           & [MWh]       \\
$P_{yakt}$                 & Production of asset $a$                                       & [MWh]       \\
$LT_{a}$                   & Lifetime of the asset $a$                                     & [years]     \\
$W_{y}^{m}$                & Weight of milestone year $y$                                  & [years]     \\
$W_{yk}^{op}$              & Weight of representative period $k$ for operational cost      & [hours]     \\
$M_{\alpha p \tau}^{op}$   & Mapping matrix between two temporal definitions $p$ and $\tau$& [-]         \\
$C^{I}$                    & Total investment cost of the energy system                    & [\EUR{}]    \\
$C^{O}$                    & Total operational cost of the energy system                   & [\EUR{}]    \\
$C_{ya}^{op}$              & Variable operational cost of asset $a$ at year $y$            & [\EUR{}/MWh]\\
$C_{ya}^{T}$               & Total investment cost of the energy asset $a$ at year $y$     & [\EUR{}/MW] \\
$SV_{ya}$                  & Salvage cost of the energy asset $a$ at year $y$              & [\EUR{}/MW] \\
$IR$                       & Interest rate                                                 & [p.u.]      \\
\hline 
\end{tabularx}

\subsection{Continuous Variables}
\begin{tabularx}{\textwidth}{| l | X | l |} 
\hline 
\textbf{Name} & \textbf{Description} & \textbf{Units}\\ 
\hline 
\endhead 
$f_{ya\alpha kt}$ & Flow from asset $a$ to asset $\alpha$                     & [MWh] \\
$r_{yarkt}$       & Reserve provision of asset $a$ to reserve requirement $r$ & [MW] \\ 
$s_{yakt}$        & Storage level of asset $a$                                & [MWh]\\
\hline 
\end{tabularx}

\subsection{Binary and Integer Variables}
The following variables are, by definition, integers; however, these definitions could be relaxed to a continuous variable for simplicity during the solution process. 
\begin{tabularx}{\textwidth}{| l | X | l |} 
\hline 
\textbf{Name} & \textbf{Description} & \textbf{Units}\\ 
\hline 
\endhead 
$u_{yakt}$           & Number of units on             & [-] \\
$\overline{u}_{ya}$  & Invested number of units       & [-] \\
\hline 
\end{tabularx}

\section{Constraints and basic methods for energy assets}
This section includes the constraints and basic methods for all energy assets $a$, i.e., producers, storage units, conversion units, and transport elements.

\subsection{Balance constraint for assets \textcolor{blue}{(One constraint to rule them all!)}}
\label{sec:balance_constraint}
This constraint is the foundation of the modelling framework in Tulipa. Equation (\ref{eq:asset_constraint}) provides a general form, but simpler versions of this constraint exist depending on the asset method. It also includes the possibility of having flexible temporal resolution on each asset \cite{Zhi2023}, allowing the model to have sectors with different temporal resolutions (e.g., hourly decisions for the electricity sector and six-hour decisions for the gas sector). In the following sections, we will explain these versions with examples to enhance comprehension.

\begin{equation}
    \label{eq:asset_constraint}
    \begin{aligned}
       \underbrace{\sum_{\alpha \in \mathcal{A}^{in}(a)}  \quad \sum_{\tau(\alpha) \in \Omega_{\alpha p \tau}} M_{\alpha p \tau}  f_{y\alpha ak\tau}}_{\textit{flows into the asset}} & \\
       + \underbrace{\sum_{\tau(a) \in \Omega_{a p \tau}} M_{a p \tau} P_{yak\tau}}_{\textit{for producer assets}}  & \\
       + \underbrace{\sum_{\tau(a) \in \Omega_{a p \tau}} M_{a p \tau} \cdot \left(s_{yak\tau} - s_{yak,\tau-1} + S_{yak\tau}^{I} \right)}_{\textit{for storage assets}}  & \\
       & \left\{\begin{array}{l} = \\ \geqslant \\ \leqslant \end{array}\right\}   \\
       & \underbrace{\sum_{\alpha \in \mathcal{A}^{out}(a)}  \quad \sum_{\tau(\alpha) \in \Omega_{\alpha p \tau}} M_{\alpha p \tau} \frac{f_{ya\alpha k\tau}}{\eta_{a \alpha}}}_{\textit{flows outgoing the asset}} \\
       & + \underbrace{\sum_{\tau(a) \in \Omega_{a p \tau}} M_{a p \tau}  D_{yak\tau}}_{\textit{for consumption assets}}  \\
       & \qquad \forall{y} \in \mathcal{Y}_m, \forall{a} \in \mathcal{A}, \forall{k} \in \mathcal{K}, \forall{p} \in \mathcal{P}
    \end{aligned}
\end{equation}

In equation (\ref{eq:asset_constraint}), the left-hand side (LHS) shows all the flows entering asset $a$ from other connected assets $\alpha$. The first term on the right-hand side (RHS) represents the outflow from the asset. This constraint has additional terms for producers, consumption, and storage assets, which will be explained later. Users can choose the constraint sense ($=$, $\leq$, or $\geq$) based on the asset's balance behavior. The constraint assumes a linear transfer function between input and output flow, which is the default method for the first version. However, we will also include an alternative method that uses a convex piece-wise linear function.

\subsubsection{Balance method for conversion assets}
In this case, equation (\ref{eq:asset_constraint}) applies to conversion assets $\mathcal{A}_{cv}$ using the balance asset method ($\mathcal{A}_{b}$). It only considers the inflows and outflows of the asset, and uses efficiency ($\eta$) as a conversion factor to calculate the outgoing flow for each incoming flow.

\begin{equation*}
    \begin{aligned}
       \sum_{\alpha \in \mathcal{A}^{in}(a)}  \quad \sum_{\tau(\alpha) \in \Omega_{\alpha p \tau}} M_{\alpha p \tau}  f_{y\alpha ak\tau}  = 
       \sum_{\alpha \in \mathcal{A}^{out}(a)}  \quad \sum_{\tau(\alpha) \in \Omega_{\alpha p \tau}} M_{\alpha p \tau} \frac{f_{ya\alpha k\tau}}{\eta_{a \alpha}} \\
       \qquad \forall{y} \in \mathcal{Y}_m, \forall{a} \in \mathcal{A}_{cv} \cap \mathcal{A}_{b}, \forall{k} \in \mathcal{K}, \forall{p} \in \mathcal{P}
    \end{aligned}
\end{equation*}

Let's take a look at an example of a fuel cell $a$ for a twelve-hourly representative period ($k$); see Figure (\ref{fig:fuel_cell_ex}):

\begin{figure}
    \centering
    \includegraphics[width=0.5\linewidth]{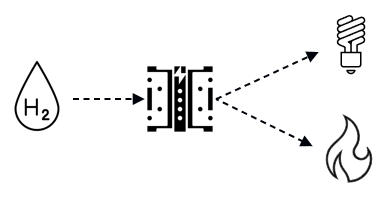}
    \caption{Conversion asset example - Fuel cell}
    \label{fig:fuel_cell_ex}
\end{figure}

\begin{itemize}
    \item Input: Hydrogen asset ($H_2$) with a resolution of 3h, i.e., $\tau(\alpha=H_2) = \{1,2,3,4\}$
    \item Output: Electricity asset ($E$) with a resolution of 1h, i.e., $\tau(\alpha=E) = \{1,2,...12\}$, and $\eta_{a,E} = 40\%$
    \item Output: Heat asset ($H$) with a resolution of 4h, i.e., $\tau(\alpha=H) = \{1,2,3\}$, and $\eta_{a,H} = 20\%$
    \item $\mathcal{A}^{in}(a)=\{H_2\}$ and $\mathcal{A}^{out}(a)=\{E,H_2\}$
\end{itemize}

The resolution determines how often we have variables for each sector (i.e., hydrogen, electricity, and heat). For instance, the hydrogen asset flow is measured by four variables ($f_{y,H_2,a,k,1}$, $f_{y,H_2,a,k,2}$, $f_{y,H_2,a,k,3}$, and $f_{y,H_2,a,k,4}$), while the total electricity generated per hour is also taken into account ($f_{y,a,E,k,1}$, $f_{y,a,E,k,2}$ ... $f_{y,a,E,k,12}$). Similarly, the fuel cell's flow to the heat asset is measured using three variables ($f_{y,a,H,k,1}$, $f_{y,a,H,k,2}$, and $f_{y,a,H,k,3}$). Please note that the resolutions are not multiples of each other and that the higher the resolution, the fewer variables. The mapping of variables to the representative period ($k$) is shown in Table \ref{tab:mapping_matrix_cv_ex}. Furthermore, we have:

\begin{table}[ht]
\centering
\begin{tabular}{|c|c|c|c|c|c|c|c|c|}
\hline
$p$ & Periods in & \textbf{$H_2$} & \textbf{$E$} & \textbf{$H$} \\ 
    & the $k$    & variables & variables & variables \\ 
\hline
\textcolor{Blue}{1} & \textcolor{Blue}{1} & \textcolor{Blue}{$1/3f_{y,H_2,a,k,1}$} & \textcolor{Blue}{$f_{y,a,E,k,1}$}  & \textcolor{Blue}{$1/4f_{y,a,H,k,1}$} \\
                    & \textcolor{Blue}{2} & \textcolor{Blue}{$1/3f_{y,H_2,a,k,1}$} & \textcolor{Blue}{$f_{y,a,E,k,2}$}  & \textcolor{Blue}{$1/4f_{y,a,H,k,1}$} \\
                    & \textcolor{Blue}{3} & \textcolor{Blue}{$1/3f_{y,H_2,a,k,1}$} & \textcolor{Blue}{$f_{y,a,E,k,3}$}  & \textcolor{Blue}{$1/4f_{y,a,H,k,1}$} \\
                    & \textcolor{Blue}{4} & \textcolor{Blue}{$1/3f_{y,H_2,a,k,2}$} & \textcolor{Blue}{$f_{y,a,E,k,4}$}  & \textcolor{Blue}{$1/4f_{y,a,H,k,1}$} \\
\hline
\textcolor{Green}{2}& \textcolor{Green}{5} & \textcolor{Green}{$1/3f_{y,H_2,a,k,2}$} & \textcolor{Green}{$f_{y,a,E,k,5}$}  & \textcolor{Green}{$1/4f_{y,a,H,k,2}$} \\
                    & \textcolor{Green}{6} & \textcolor{Green}{$1/3f_{y,H_2,a,k,2}$} & \textcolor{Green}{$f_{y,a,E,k,6}$}  & \textcolor{Green}{$1/4f_{y,a,H,k,2}$} \\
                    & \textcolor{Green}{7} & \textcolor{Green}{$1/3f_{y,H_2,a,k,3}$} & \textcolor{Green}{$f_{y,a,E,k,7}$}  & \textcolor{Green}{$1/4f_{y,a,H,k,2}$} \\
                    & \textcolor{Green}{8} & \textcolor{Green}{$1/3f_{y,H_2,a,k,3}$} & \textcolor{Green}{$f_{y,a,E,k,8}$}  & \textcolor{Green}{$1/4f_{y,a,H,k,2}$} \\
\hline
\textcolor{Plum}{3} & \textcolor{Plum}{9 } & \textcolor{Plum}{$1/3f_{y,H_2,a,k,3}$} & \textcolor{Plum}{$f_{y,a,E,k,9}$ } & \textcolor{Plum}{$1/4f_{y,a,H,k,3}$} \\
                    & \textcolor{Plum}{10} & \textcolor{Plum}{$1/3f_{y,H_2,a,k,4}$} & \textcolor{Plum}{$f_{y,a,E,k,10}$} & \textcolor{Plum}{$1/4f_{y,a,H,k,3}$} \\
                    & \textcolor{Plum}{11} & \textcolor{Plum}{$1/3f_{y,H_2,a,k,4}$} & \textcolor{Plum}{$f_{y,a,E,k,11}$} & \textcolor{Plum}{$1/4f_{y,a,H,k,3}$} \\
                    & \textcolor{Plum}{12} & \textcolor{Plum}{$1/3f_{y,H_2,a,k,4}$} & \textcolor{Plum}{$f_{y,a,E,k,12}$} & \textcolor{Plum}{$1/4f_{y,a,H,k,3}$} \\
\hline
\end{tabular}
\caption{Mapping between variables and periods for the fuel-cell $a$ example}
\label{tab:mapping_matrix_cv_ex}
\end{table}

\begin{itemize}
    \item The index $p$ represents the maximum resolution of connections to the asset. In this case, $p = \{1,2,3\}$ and each value of $p$ will span 4h. Table \ref{tab:mapping_matrix_cv_ex} uses three different colours to highlight the periods in the $k$ that are considered for each value of $p$.
    \item The matrices $M_{\alpha p \tau}$ can be obtained from Table \ref{tab:mapping_matrix_cv_ex}. These matrices will map the values from the auxiliary index $p$ to each time resolution.
\end{itemize}

\[
M_{H_2,p,\tau(H_2)}= \begin{array}{l} \textcolor{lightgray}{p=1} \\ \textcolor{lightgray}{p=2} \\ \textcolor{lightgray}{p=3} \end{array}
\left(\begin{array}{cccc}
\textcolor{Blue}{1} & \textcolor{Blue}{1/3} & 0   & 0 \\
0 & \textcolor{Green}{2/3} & \textcolor{Green}{2/3} & 0 \\
0 & 0   & \textcolor{Plum}{1/3} & \textcolor{Plum}{1} \\
\end{array}\right)
\]

\[
M_{E,p,\tau(E)}=\begin{array}{l} \textcolor{lightgray}{p=1} \\ \textcolor{lightgray}{p=2} \\ \textcolor{lightgray}{p=3} \end{array}
\left(\begin{array}{cccccccccccc}
\textcolor{Blue}{1} & \textcolor{Blue}{1} & \textcolor{Blue}{1} & \textcolor{Blue}{1} & 0 & 0 & 0 & 0 & 0 & 0 & 0 & 0\\
0 & 0 & 0 & 0 & \textcolor{Green}{1} & \textcolor{Green}{1} & \textcolor{Green}{1} & \textcolor{Green}{1} & 0 & 0 & 0 & 0\\
0 & 0 & 0 & 0 & 0 & 0 & 0 & 0 & \textcolor{Plum}{1} & \textcolor{Plum}{1} & \textcolor{Plum}{1} & \textcolor{Plum}{1}\\
\end{array}\right)
\]

\[
M_{H,p,\tau(H)}=\begin{array}{l} \textcolor{lightgray}{p=1} \\ \textcolor{lightgray}{p=2} \\ \textcolor{lightgray}{p=3} \end{array}
\left(\begin{array}{ccc}
\textcolor{Blue}{1} & 0 & 0 \\
0 & \textcolor{Green}{1} & 0 \\
0 & 0 & \textcolor{Plum}{1} \\
\end{array}\right)
\]

The coefficients in the matrices determine the fraction of each variable taken into account for each constraint. For example, when $p=1$, the complete value of $f_{y,H_2,a,k,1}$ is considered, but only $1/3$ of $f_{y,H_2,a,k,2}$ is taken into account since $p$ represents a four-hour interval (refer to Table \ref{tab:mapping_matrix_cv_ex}). Similarly, when $p=2$, $2/3$ of $f_{y,H_2,k,2}$ and $2/3$ of $f_{y,H_2,a,k,3}$ are considered, and so on.

$p=1$:
\begin{equation*}
f_{y,H_2,a,k,1}+\frac{1}{3}f_{y,H_2,a,k,2} = \frac{f_{y,a,E,k,1}}{0.4} + \frac{f_{y,a,E,k,2}}{0.4} + \frac{f_{y,a,E,k,3}}{0.4} + \frac{f_{y,a,E,k,4}}{0.4} + \frac{f_{y,a,H,k,1}}{0.2}
\end{equation*}

$p=2$:
\begin{equation*}
\frac{2}{3}f_{y,H_2,a,k,2}+\frac{2}{3}f_{y,H_2,a,k,3} = \frac{f_{y,a,E,k,5}}{0.4} + \frac{f_{y,a,E,k,6}}{0.4} + \frac{f_{y,a,E,k,7}}{0.4} + \frac{f_{y,a,E,k,8}}{0.4} + \frac{f_{y,a,H,k,2}}{0.2}
\end{equation*}

$p=3$:
\begin{equation*}
\frac{1}{3}f_{y,H_3,a,k,1}+f_{y,H_2,a,k,4} = \frac{f_{y,a,E,k,9}}{0.4} + \frac{f_{y,a,E,k,10}}{0.4} + \frac{f_{y,a,E,k,11}}{0.4} + \frac{f_{y,a,E,k,12}}{0.4} + \frac{f_{y,a,H,k,3}}{0.2}
\end{equation*}

To define the sets $\Omega_{\alpha p \tau}$, we use a tuple of three elements derived from the matrices $M_{\alpha p \tau}$ so that the positive values compound the elements in the set.

\begin{equation*}
\Omega_{\alpha p \tau} = \{ (\alpha,p,\tau) | M_{\alpha p \tau} > 0 \}
\end{equation*}

\subsubsection{Balance method for production assets}
Production assets are elements that produce energy from a resource and distribute it to one or more energy assets in the network. These assets can include solar, wind, nuclear, and energy imports. When the method is applied to production assets, equation (\ref{eq:asset_constraint}) can be simplified into the following expression.

\begin{equation*}
    \begin{aligned}
       \sum_{\tau(a) \in \Omega_{a p \tau}} M_{a p \tau} P_{yak\tau}  
       \left\{\begin{array}{l} = \\ \geqslant \\ \leqslant \end{array}\right\}   
       \sum_{\alpha \in \mathcal{A}^{out}(a)}  \quad \sum_{\tau(\alpha) \in \Omega_{\alpha p \tau}} M_{\alpha p \tau} f_{ya\alpha k\tau} \\
       \qquad \forall{y} \in \mathcal{Y}_m, \forall{a} \in \mathcal{A}_{p} \cap \mathcal{A}_{b}, \forall{k} \in \mathcal{K}, \forall{p} \in \mathcal{P}
    \end{aligned}
\end{equation*}

\begin{figure}
    \centering
    \includegraphics[width=0.25\linewidth]{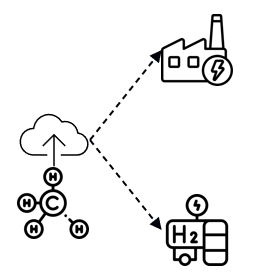}
    \caption{Production asset example - Methane producer}
    \label{fig:methane_producer_ex}
\end{figure}

On the LHS, we have the representation of the production resource availability, while on the RHS, we have the flows that are directed towards other assets. It is important to note that the efficiency is equal to one ($\eta_{\alpha a}=1$) due to the absence of any energy conversion in the production assets. This approach can prove helpful for various production assets, but it's particularly relevant for energy imports. As an example, take a look at Figure \ref{fig:methane_producer_ex}, which showcases a methane source that can provide fuel to a gas-fired plant and a $H_2$ generator like a steam-methane reformation (SMR), having a twelve-hourly representative period ($k$) again:

\begin{itemize}
    \item Production: Methane ($CH_4$) with a resolution of 6h, i.e., $\tau(\alpha=CH_4) = \{1,2\}$
    \item Output: Gas turbine ($GT$) with a resolution of 3h, i.e., $\tau(\alpha=GT) = \{1,2,3,4\}$
    \item Output: Steam-methane reformation ($SMR$) with a resolution of 4h, i.e., $\tau(\alpha=SMR) = \{1,2,3\}$
    \item $\mathcal{A}^{out}(CH_4)=\{GT,SMR\}$
\end{itemize}

As shown in the previous section, we can map variables to the representative period ($k$) using Table \ref{tab:mapping_matrix_p_ex}. The index $p$ represents the maximum resolution of connections to the asset, with $p$ being either 1 or 2, each representing a 6-hour span. Table \ref{tab:mapping_matrix_p_ex} uses two colours to highlight the periods in the $p$ considered for each $p$ value. The $M_{\alpha p \tau}$ matrices can be derived from the table values.

\begin{table}[ht]
\centering
\begin{tabular}{|c|c|c|c|c|c|c|c|c|}
\hline
$p$ & Periods in & \textbf{$CH_4$} & \textbf{$GT$} & \textbf{$SMR$} \\ 
    & the $k$   &  profile values & variables     & variables \\ 
\hline
\textcolor{Blue}{1} & \textcolor{Blue}{1} & \textcolor{Blue}{$1/6P_{y,CH_4,k,1}$} & \textcolor{Blue}{$1/3f_{y,CH_4,GT,k,1}$}  & \textcolor{Blue}{$1/4f_{y,CH_4,SMR,k,1}$} \\
                    & \textcolor{Blue}{2} & \textcolor{Blue}{$1/6P_{y,CH_4,k,1}$} & \textcolor{Blue}{$1/3f_{y,CH_4,GT,k,1}$}  & \textcolor{Blue}{$1/4f_{y,CH_4,SMR,k,1}$} \\
                    & \textcolor{Blue}{3} & \textcolor{Blue}{$1/6P_{y,CH_4,k,1}$} & \textcolor{Blue}{$1/3f_{y,CH_4,GT,k,1}$}  & \textcolor{Blue}{$1/4f_{y,CH_4,SMR,k,1}$} \\
                    & \textcolor{Blue}{4} & \textcolor{Blue}{$1/6P_{y,CH_4,k,1}$} & \textcolor{Blue}{$1/3f_{y,CH_4,GT,k,2}$}  & \textcolor{Blue}{$1/4f_{y,CH_4,SMR,k,1}$} \\
                    & \textcolor{Blue}{5} & \textcolor{Blue}{$1/6P_{y,CH_4,k,1}$} & \textcolor{Blue}{$1/3f_{y,CH_4,GT,k,2}$}  & \textcolor{Blue}{$1/4f_{y,CH_4,SMR,k,2}$} \\
                    & \textcolor{Blue}{6} & \textcolor{Blue}{$1/6P_{y,CH_4,k,1}$} & \textcolor{Blue}{$1/3f_{y,CH_4,GT,k,2}$}  & \textcolor{Blue}{$1/4f_{y,CH_4,SMR,k,2}$} \\
\hline
\textcolor{Green}{2}& \textcolor{Green}{7} & \textcolor{Green}{$1/6P_{y,CH_4,k,2}$} & \textcolor{Green}{$1/3f_{y,CH_4,GT,k,3}$}  & \textcolor{Green}{$1/4f_{y,CH_4,SMR,k,2}$} \\
                    & \textcolor{Green}{8} & \textcolor{Green}{$1/6P_{y,CH_4,k,2}$} & \textcolor{Green}{$1/3f_{y,CH_4,GT,k,3}$}  & \textcolor{Green}{$1/4f_{y,CH_4,SMR,k,2}$} \\
                    & \textcolor{Green}{9} & \textcolor{Green}{$1/6P_{y,CH_4,k,2}$} & \textcolor{Green}{$1/3f_{y,CH_4,GT,k,3}$}  & \textcolor{Green}{$1/4f_{y,CH_4,SMR,k,3}$} \\
                    & \textcolor{Green}{10}& \textcolor{Green}{$1/6P_{y,CH_4,k,2}$} & \textcolor{Green}{$1/3f_{y,CH_4,GT,k,4}$}  & \textcolor{Green}{$1/4f_{y,CH_4,SMR,k,3}$} \\
                    & \textcolor{Green}{11}& \textcolor{Green}{$1/6P_{y,CH_4,k,2}$} & \textcolor{Green}{$1/3f_{y,CH_4,GT,k,4}$}  & \textcolor{Green}{$1/4f_{y,CH_4,SMR,k,3}$} \\
                    & \textcolor{Green}{12}& \textcolor{Green}{$1/6P_{y,CH_4,k,2}$} & \textcolor{Green}{$1/3f_{y,CH_4,GT,k,4}$}  & \textcolor{Green}{$1/4f_{y,CH_4,SMR,k,3}$} \\
\hline
\end{tabular}
\caption{Mapping between variables and periods for the methane imports $CH_4$ example}
\label{tab:mapping_matrix_p_ex}
\end{table}

\[
M_{CH_4,k,\tau(CH_4)}= \begin{array}{l} \textcolor{lightgray}{p=1} \\ \textcolor{lightgray}{p=2} \end{array}
\left(\begin{array}{cc}
\textcolor{Blue}{1} & 0  \\
0 & \textcolor{Green}{1} \\
\end{array}\right)
\]

\[
M_{GT,k,\tau(GT)}=\begin{array}{l} \textcolor{lightgray}{p=1} \\ \textcolor{lightgray}{p=2} \end{array}
\left(\begin{array}{cccc}
\textcolor{Blue}{1} & \textcolor{Blue}{1} & 0 & 0\\
0 & 0 & \textcolor{Green}{1} & \textcolor{Green}{1} \\
\end{array}\right)
\]

\[
M_{SMR,k,\tau(SMR)}=\begin{array}{l} \textcolor{lightgray}{p=1} \\ \textcolor{lightgray}{p=2} \end{array}
\left(\begin{array}{ccc}
\textcolor{Blue}{1} & \textcolor{Blue}{1/2} & 0 \\
0 & \textcolor{Green}{1/2} & \textcolor{Green}{1} \\
\end{array}\right)
\]

In this example, the equation (\ref{eq:asset_constraint}) can be used with different options to indicate minimum production, maximum production, or fixed production. We will use the maximum production limit as an example in this case.

$p=6$:
\begin{equation*}
    \begin{aligned}
    P_{y,CH_4,k,1} \geq & f_{y,CH_4,GT,k,1} + f_{y,CH_4,GT,k,2} + f_{y,CH_4,SMR,k,1} + \frac{1}{2} f_{y,CH_4,SMR,k,2}
    \end{aligned}
\end{equation*}

$p=12$:
\begin{equation*}
    \begin{aligned}
    P_{y,CH_4,k,2} \geq & f_{y,CH_4,GT,k,3} + f_{y,CH_4,GT,k,4} + \frac{1}{2} f_{y,CH_4,SMR,k,2} + f_{y,CH_4,SMR,k,3}
    \end{aligned}
\end{equation*}

This method is beneficial for representing production assets that are linked to multiple assets with varying time resolutions, as illustrated in the example. However, all production assets don't need to utilize this method. There are other specific methods to restrict the flows from the production asset to the grid if there is only one connection and the production assets can offer reserves, as discussed in Section \ref{sec:operation_limits}.

\subsubsection{Balance method for consumption assets}
Consumption assets are elements that demand energy from one or more energy assets in the network. These assets can include energy demands and exports. When the method is applied to consumption assets, equation (\ref{eq:asset_constraint}) can be simplified into the following expression.

\begin{equation*}
    \begin{aligned}
       \sum_{\alpha \in \mathcal{A}^{in}(a)}  \quad \sum_{\tau(\alpha) \in \Omega_{\alpha p \tau}} M_{\alpha p \tau}  f_{y\alpha ak\tau} \left\{\begin{array}{l} = \\ \geqslant \\ \leqslant \end{array}\right\}
       \sum_{\tau(a) \in \Omega_{a p \tau}} M_{a p \tau}  D_{yak\tau}  \\
       \qquad \forall{y} \in \mathcal{Y}_m, \forall{a} \in \mathcal{A}_{c} \cap \mathcal{A}_{b}, \forall{k} \in \mathcal{K}, \forall{p} \in \mathcal{P}
    \end{aligned}
\end{equation*}

\begin{figure}
    \centering
    \includegraphics[width=0.25\linewidth]{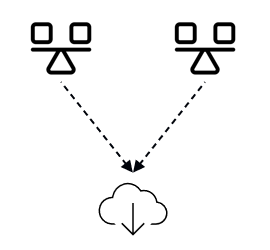}
    \caption{Consumption asset example - Electricity exports}
    \label{fig:exports_consumption_ex}
\end{figure}

On the LHS, we have the representation of the flows that are directed towards the consumption asset, while on the RHS, we have the consumption demand profile. This approach can prove helpful for various consumption assets, but it's particularly relevant for energy exports. As an example, take a look at Figure \ref{fig:exports_consumption_ex}, which showcases electricity exports that can take their demand from two nodes/countries: 

\begin{itemize}
    \item Input: Node 1 ($N1$) with a resolution of 1h, i.e., $\tau(\alpha=N1) = \{1,2...12\}$
    \item Input: Node 2 ($N2$) with a resolution of 4h, i.e., $\tau(\alpha=N2) = \{1,2,3\}$
    \item Consumption: Electricity exports ($E$) with a resolution of 6h, i.e., $\tau(\alpha=E) = \{1,2\}$
    \item $\mathcal{A}^{in}(E)=\{N1,N2\}$
\end{itemize}

Table \ref{tab:mapping_matrix_p_ex} maps the variables to the periods in the representative period ($k$). The $M_{\alpha p \tau}$ matrices can be derived from the table values, such as before.

\begin{table}[ht]
\centering
\begin{tabular}{|c|c|c|c|c|c|c|c|c|}
\hline
$p$ & Periods in & \textbf{$N1$} & \textbf{$N2$} & \textbf{$D$} \\ 
    & the $k$   & variables & variables &  profile values \\ 
\hline
\textcolor{Blue}{1} & \textcolor{Blue}{1}  & \textcolor{Blue}{$f_{y,N1,E,k,1}$}  & \textcolor{Blue}{$1/4f_{y,N2,E,k,1}$} & \textcolor{Blue}{$1/6D_{y,E,k,1}$} \\
                    & \textcolor{Blue}{2}  & \textcolor{Blue}{$f_{y,N1,E,k,2}$}  & \textcolor{Blue}{$1/4f_{y,N2,E,k,1}$} & \textcolor{Blue}{$1/6D_{y,E,k,1}$} \\
                    & \textcolor{Blue}{3}  & \textcolor{Blue}{$f_{y,N1,E,k,3}$}  & \textcolor{Blue}{$1/4f_{y,N2,E,k,1}$} & \textcolor{Blue}{$1/6D_{y,E,k,1}$} \\
                    & \textcolor{Blue}{4}  & \textcolor{Blue}{$f_{y,N1,E,k,4}$}  & \textcolor{Blue}{$1/4f_{y,N2,E,k,1}$} & \textcolor{Blue}{$1/6D_{y,E,k,1}$} \\
                    & \textcolor{Blue}{5}  & \textcolor{Blue}{$f_{y,N1,E,k,5}$}  & \textcolor{Blue}{$1/4f_{y,N2,E,k,2}$} & \textcolor{Blue}{$1/6D_{y,E,k,1}$} \\
                    & \textcolor{Blue}{6}  & \textcolor{Blue}{$f_{y,N1,E,k,6}$}  & \textcolor{Blue}{$1/4f_{y,N2,E,k,2}$} & \textcolor{Blue}{$1/6D_{y,E,k,1}$} \\
\hline
\textcolor{Green}{2}& \textcolor{Green}{7} & \textcolor{Green}{$f_{y,N1,E,k,7}$} & \textcolor{Green}{$1/4f_{y,N2,E,k,2}$}& \textcolor{Green}{$1/6D_{y,E,k,2}$} \\
                    & \textcolor{Green}{8} & \textcolor{Green}{$f_{y,N1,E,k,8}$} & \textcolor{Green}{$1/4f_{y,N2,E,k,2}$}& \textcolor{Green}{$1/6D_{y,E,k,2}$} \\
                    & \textcolor{Green}{9} & \textcolor{Green}{$f_{y,N1,E,k,9}$} & \textcolor{Green}{$1/4f_{y,N2,E,k,3}$}& \textcolor{Green}{$1/6D_{y,E,k,2}$} \\
                    & \textcolor{Green}{10}& \textcolor{Green}{$f_{y,N1,E,k,10}$}& \textcolor{Green}{$1/4f_{y,N2,E,k,3}$}& \textcolor{Green}{$1/6D_{y,E,k,2}$} \\
                    & \textcolor{Green}{11}& \textcolor{Green}{$f_{y,N1,E,k,11}$}& \textcolor{Green}{$1/4f_{y,N2,E,k,3}$}& \textcolor{Green}{$1/6D_{y,E,k,2}$} \\
                    & \textcolor{Green}{12}& \textcolor{Green}{$f_{y,N1,E,k,12}$}& \textcolor{Green}{$1/4f_{y,N2,E,k,3}$}& \textcolor{Green}{$1/6D_{y,E,k,2}$} \\
\hline
\end{tabular}
\caption{Mapping between variables and periods for the electricity exports $E$ example}
\label{tab:mapping_matrix_c_ex}
\end{table}

\[
M_{N1,k,\tau(N1)}= \begin{array}{l} \textcolor{lightgray}{p=1} \\ \textcolor{lightgray}{p=2} \end{array}
\left(\begin{array}{cccccccccccc}
\textcolor{Blue}{1} & \textcolor{Blue}{1} & \textcolor{Blue}{1} & \textcolor{Blue}{1} & \textcolor{Blue}{1} & \textcolor{Blue}{1} &  0 & 0 & 0 & 0 & 0 & 0 \\
0 & 0 & 0 & 0 & 0 & 0 & \textcolor{Green}{1} & \textcolor{Green}{1} & \textcolor{Green}{1} & \textcolor{Green}{1} & \textcolor{Green}{1} & \textcolor{Green}{1} \\
\end{array}\right)
\]

\[
M_{N2,k,\tau(N2)}=\begin{array}{l} \textcolor{lightgray}{p=1} \\ \textcolor{lightgray}{p=2} \end{array}
\left(\begin{array}{ccc}
\textcolor{Blue}{1} & \textcolor{Blue}{1/2} & 0 \\
0 & \textcolor{Green}{1/2} & \textcolor{Green}{1} \\
\end{array}\right)
\]

\[
M_{D,k,\tau(D)}=\begin{array}{l} \textcolor{lightgray}{p=1} \\ \textcolor{lightgray}{p=2} \end{array}
\left(\begin{array}{cc}
\textcolor{Blue}{1} & 0  \\
0 & \textcolor{Green}{1} \\
\end{array}\right)
\]

In this example, the equation (\ref{eq:asset_constraint}) can be used with different options to indicate minimum consumption, maximum consumption, or fixed consumption. We will use the minimum consumption limit as an example in this case.

$p=6$:
\begin{equation*}
    \begin{aligned}
    &f_{y,N1,E,k,1} + f_{y,N1,E,k,2} + f_{y,N1,E,k,3} + f_{y,N1,E,k,4}  \\
    &    + f_{y,N1,E,k,5} + f_{y,N1,E,k,6} + f_{y,N2,E,k,1} + \frac{1}{2} f_{y,N2,E,k,2}  \leq D_{y,E,k,1} \\
    \end{aligned}
\end{equation*}

$p=12$:
\begin{equation*}
    \begin{aligned}
    &f_{y,N1,E,k,7} + f_{y,N1,E,k,8} + f_{y,N1,E,k,9} + f_{y,N1,E,k,10}  \\
    &    + f_{y,N1,E,k,11} + f_{y,N1,E,k,12} + \frac{1}{2} f_{y,N2,E,k,2} + f_{y,N2,E,k,3}  \leq D_{y,E,k,2} \\
    \end{aligned}
\end{equation*}

This approach helps depict consumption assets that are associated with several assets having different time resolutions, as shown in the example. Rather than using a node to combine the flows, we can simplify the connections by directing all the flows towards the consumption asset that has a demand, which ultimately reduces the number of variables.

\subsubsection{Balance method for transport assets}
The RHS of the equation represents all the production that goes into the transport asset. The first term on the LHS is the flow going out of the transport asset. Since transport assets do not convert energy, the efficiency (represented by the symbol $\eta$) equals 1. Additionally, the constraint is an equality constraint because there is no production or consumption of energy at the transport assets. We can represent, for example, a nodal balance in an intermediate point, such as in a transhipment problem.

\begin{equation*}
    \begin{aligned}
       \sum_{\alpha \in \mathcal{A}^{in}(a)}   \quad \sum_{\tau(\alpha) \in \Omega_{\alpha p \tau}} M_{\alpha p \tau}  f_{y\alpha ak\tau} =
       \sum_{\alpha \in \mathcal{A}^{out}(a)}  \quad \sum_{\tau(\alpha) \in \Omega_{\alpha p \tau}} M_{\alpha p \tau} \frac{f_{ya\alpha k\tau}}{\eta_{a \alpha}} \\
       \qquad \forall{y} \in \mathcal{Y}_m, \forall{a} \in \mathcal{A}_{t} \cap \mathcal{A}_{b}, \forall{k} \in \mathcal{K}, \forall{p} \in \mathcal{P}
    \end{aligned}
\end{equation*}

\subsubsection{Balance method for storage assets}
For storage assets using a balance method, the LHS represents the charging flows into the storage asset, typically only one. The LHS also has the change in the storage level of the asset and inflows (if applicable). On the RHS, there are the discharging flows (also typically only one). Each equation component may have a different time resolution, allowing for situations where the discharging and charging decision variables have a lower resolution than the desired storage level variables. For example, hourly charging/discharging decisions with a 24-hour resolution for the storage level help represent pumped-hydro storage units.

\begin{equation*}
    \begin{aligned}
       \sum_{\alpha \in \mathcal{A}^{in}(a)}  \quad \sum_{\tau(\alpha) \in \Omega_{\alpha p \tau}} M_{\alpha p \tau}  f_{y\alpha ak\tau} 
       & + \sum_{\tau(a) \in \Omega_{a p \tau}} M_{a p \tau} \cdot \left(s_{yak\tau} - s_{yak,\tau-1} + S_{yak\tau}^{I} \right)  
       = \\
       & \sum_{\alpha \in \mathcal{A}^{out}(a)} \quad \sum_{\tau(\alpha) \in \Omega_{\alpha p \tau}} M_{\alpha p \tau} \frac{f_{ya\alpha k\tau}}{\eta_{a \alpha}} \\
       & \qquad \forall{y} \in \mathcal{Y}_m, \forall{a} \in \mathcal{A}_{s} \cap \mathcal{A}_{b}, \forall{k} \in \mathcal{K}, \forall{p} \in \mathcal{P}
    \end{aligned}
\end{equation*}

%%See Section \ref{sec:inter-period-storage} for the inter-storage balance.

\subsection{Operation limit constraints}
\label{sec:operation_limits}
Constraints (\ref{eq:max_prod}) and (\ref{eq:min_prod}) represent the maximum and minimum output flow for each asset. Note that some assets may have only one output flow (e.g. renewable assets), while others may have several (e.g., fuel cells). Unlike the balance constraint in Section \ref{sec:balance_constraint}, these constraints allow for the consideration of reserve ($r$) and unit commitment ($u$) variables. Furthermore, the parameters of the profile, which are represented by $\tilde{F}^{max}$ and $\tilde{F}^{max}$, can be set with a temporal resolution that is separate from that of the variables. For example, the profiles could have hourly values while the variable definition is set every 4 hours. To account for this, the formulation includes a mapping matrix that multiplies the profiles, as seen in the balance constraint.

\begin{equation}
 \label{eq:max_prod}
 \begin{aligned}
     \sum_{\alpha \in \mathcal{A}^{out}(a)}f_{ya\alpha k t}
    +\sum_{r \in \mathcal{R}_{a}^{+}} r_{yarkt}  \leq &
    \overline{F}_{a} u_{yakt} \sum_{p \in \Omega_{a p t}} M_{a p t} \tilde{F}_{yakp}^{max}
    \\ & \forall{y} \in \mathcal{Y}_m, \forall{a} \in \mathcal{A}, \forall{k} \in \mathcal{K}, \forall{t} \in \tau(a)
 \end{aligned}
\end{equation}

\begin{equation}
 \label{eq:min_prod}
  \begin{aligned}
     \sum_{\alpha \in \mathcal{A}^{out}(a)}f_{ya\alpha k t}
    -\sum_{r \in \mathcal{R}_{a}^{-}} r_{yarkt} \geq &
    \overline{F}_{a} u_{yakt} \sum_{p \in \Omega_{a p t}} M_{a p t} \tilde{F}_{yakp}^{min} 
    \\ &  \forall{y} \in \mathcal{Y}_m, \forall{a} \in \mathcal{A}, \forall{k} \in \mathcal{K}, \forall{t} \in \tau(a)
 \end{aligned}
\end{equation}

As a general feature, the minimum output flow constraint helps represent minimum production levels without a unit commitment method, e.g., a nuclear plant with a minimum stable load of 80\%; however, it can go up to 100\% of its capacity if needed. Notice that for assets with a unit commitment method, then $\tilde{F}^{min}=0$, and for assets with a transport method, then $\tilde{F}^{min}=-1$. The upwards/downwards spinning reserve variables will appear or not in the constraint depending on whether the asset has the reserve method assigned or not (e.g., transport assets do not have reserve variables, whereas power plant units might have the possibility to provide different types of reserves).

\subsection{Reserve requirement constraints}
\label{sec:reserve_requirement}
The reserve limit ensures that all assets capable of providing reserve type $r$, denoted by $\Omega_{ar}$, must sum up to at least the minimum reserve requirement ($R$) at each time period ($p$) in the representative period ($k$). Note that the constraint is defined on the temporal definition of the parameter $R$ (i.e., $p$). Therefore, we need to multiply the variable for reserve with the mapping matrix between the temporal definitions of the variables and the parameter.  

\begin{equation}
 \label{eq:reserve_requirement}
    \sum_{a \in \Omega_{ar}} \quad \sum_{\tau(\alpha) \in \Omega_{\alpha p \tau}} M_{\alpha p \tau} r_{yark\tau} \geq R_{yrkp}
    \quad  \forall{y} \in \mathcal{Y}_m, \forall{r} \in \mathcal{R}, \forall{k} \in \mathcal{K}, \forall{p} \in \mathcal{P}
\end{equation}

\subsection{Investment methods}
The subset $\mathcal{A}_{i}$ includes all assets with an investment method. Therefore, assets can have two types of constraints: with and without investment.

\subsubsection{Asset with investment}
The following constraints relate the available energy asset units with the investment decision variables and limit the maximum investment. For instance, constraint (\ref{eq:asset_with_investment}) imposes that the available energy asset units at year $y$ must be lower or equal to those available at the beginning plus those built in the previous years. In addition, constraint (\ref{eq:investment_potential}) ensure that the cumulative investment is lower or equal to the investment potential $\overline{F}_{ya}^{max}$ at a specific year $y$ for each energy asset.

\begin{equation}
 \label{eq:asset_with_investment}
    u_{yakt} \leq \overline{U}_{ya} + \sum_{\psi=\max(y-LT_a+1,0)}^{y}\overline{u}_{\psi,a}  \quad \forall{y} \in \mathcal{Y}_m, \forall{a} \in \mathcal{A}_{i}, \forall{k} \in \mathcal{K}, \forall{t}  \in \tau(a)
\end{equation}

\begin{equation}
 \label{eq:investment_potential}
    \overline{F}_{a} \sum_{\psi=\max(y-LT_a+1,0)}^{y}\overline{u}_{\psi a} \leq  \overline{F}_{ya}^{max}  \quad \forall{y} \in \mathcal{Y}_m, \forall{a} \in \mathcal{A}_{i}
\end{equation}

\subsubsection{Asset without investment}
For assets excluded from the investment method ($a \notin \mathcal{A}_{i}$), the available units are only limited by their initial amount and $\overline{u}$ is set to 0.

\begin{equation}
 \label{eq:asset_without_investment}
    u_{yakt} \leq \overline{U}_{ya} \quad \forall{y} \in \mathcal{Y}_m, \forall{a} \notin \mathcal{A}_{i}, \forall{k} \in \mathcal{K}, \forall{t}  \in \tau(a)
\end{equation}

\section{Objective function}
The objective function is to minimise the energy system's total investment and operational cost.

\begin{equation}
 \label{eq:ob_total_cost}
  \begin{aligned}
    & \min \quad C^{I} + C^{O}  \\
    & C^{I} = \sum_{y \in \mathcal{Y}_m}\sum_{a \in \mathcal{A}_{i}}\frac{1}{(1+IR)^{y}} \cdot (C_{ya}^{T}-SV_{ya})\cdot \overline{F}_{a} \overline{u}_{ya} \\ 
    & C^{O} = \sum_{y \in \mathcal{Y}_m}\frac{1}{(1+IR)^{y}} \cdot W_{y}^{m} \cdot \sum_{a \in \mathcal{A}} C_{ya}^{op} \cdot \sum_{\alpha \in \mathcal{A}^{out}(a)}\sum_{k \in \mathcal{K}}W_{yk}^{op} \cdot \sum_{t \in \mathcal{T}} \sum_{\tau(\alpha) \in \Omega_{\alpha t \tau}} M_{\alpha t \tau} f_{ya\alpha k\tau}
  \end{aligned}
\end{equation} 

The term $C^{I}$ represents the investment costs of energy assets. It is multiplied by the interest rate $IR$ for each milestone year ($y \in \mathcal{Y}_m$). It considers the total investment cost $C_y^T$ and salvage value $SV_y$ for the end-of-horizon effect.

The term $C^{O}$ represents the total variable production cost $C_{ya}^{op}$, also known as operational cost. The model will include other costs like start-up costs, no-load costs, and shut-down costs once unit commitment methods are defined. Operational cost includes the interest rate $IR$, the weight of the milestone year $W_{y}^{m}$, and the weight of the representative period $W_{yk}^{op}$. Each flow in the model is multiplied by the mapping matrix to account for the variable decision's suitable duration due to the temporal structure's flexibility.

The article referenced by \cite{Tejada2023} explores various methods for modeling multiyear investments in energy systems, outlining the pros and cons of each one.

\section{GitHub repository}
The Tulipa Energy Model is implemented in Julia \cite{Julia-2017} using the JuMP \cite{Lubin2023} and Graphs \cite{Graphs2021} packages. The code is available on GitHub and is based on the formulation presented in this document.
\\
\\
GitHub repository: \href{https://github.com/TulipaEnergy/TulipaEnergyModel.jl}{https://github.com/TulipaEnergy/TulipaEnergyModel.jl}

\printbibliography

\end{document}